\documentclass{amsart}
\usepackage{amsmath}
\usepackage{amssymb}

\newenvironment{preuve}{\noindent {\it Proof}}
{\hfill $\square$ \vspace{0.25cm}}

\newcommand{\ds}{\displaystyle}
\newcommand{\rr}{{\mathbb{R}}}
\newcommand{\E}{{\mathbb{E}}}
\newcommand{\sg}{{\rm{sg}}}
\newcommand{\intot}{\ds\int_0^t}
\newcommand{\indiq}{{\bf 1}}
\newcommand{\lb}{\left< }
\newcommand{\rb}{\right> }
\newcommand{\ala}{\nonumber \\}
\newcommand{\cL}{{{\mathcal L}}}
\newcommand{\cX}{{{\mathcal X}}}
\newcommand{\cF}{{{\mathcal F}}}
\newcommand{\cH}{{{\mathcal H}}}
\newcommand{\cI}{{{\mathcal I}}}
\newcommand{\cB}{{{\mathcal B}}}
\newcommand{\e}{\epsilon}
\newcommand{\vip}{\vskip0.2cm}
\newcommand{\tQ}{{\tilde Q}}
\newcommand{\tu}{{\tilde u}}

\newtheorem{theo}{Theorem}
\newtheorem{pro}[theo]{Proposition}

\newtheorem{lem}[theo]{Lemma}
\newtheorem{defi}[theo]{Definition}
\newtheorem{cor}[theo]{Corollary}

\begin{document}

\title{Stability of the stochastic heat equation in $L^1([0,1])$}
\author{Nicolas Fournier}
\author{Jacques Printems}
\address{Universit\'e Paris Est, LAMA UMR 8050, 
Facult\'e des Sciences et Technologies, 61 avenue du G\'en\'eral de Gaulle, 
94010 Cr\'eteil Cedex, France.}
\email{nicolas.fournier@univ-paris12.fr, printems@univ-paris12.fr}

\def\abstractname{Abstract}
\begin{abstract}
\noindent We consider the 
white-noise driven stochastic heat equation on $[0,\infty)\times[0,1]$
with Lipschitz-continuous
drift and diffusion coefficients $b$ and $\sigma$.
We derive an inequality for the $L^1([0,1])$-norm of the difference
between two solutions.
Using some martingale arguments, we show that this inequality provides some 
{\it a priori} estimates on solutions. This allows us to prove
the strong existence and (partial) uniqueness
of weak solutions when the initial condition 
belongs only to $L^1([0,1])$, and the stability of the solution 
with respect to this initial condition. 
We also obtain, under some conditions, 
some results concerning the large time behavior
of solutions: uniqueness of the possible 
invariant distribution and asymptotic confluence of solutions.
\end{abstract}



\maketitle

\section{Introduction and results}\label{intro}

\subsection{The equation}
Consider the stochastic heat equation with Neumann boundary 
conditions:
\begin{equation}\label{she}
\left\{
\begin{array}{rcl}
\ds\partial_t u(t,x) &=& 
\partial_{xx}u(t,x) + b(u(t,x)) + \sigma(u(t,x))\dot W(t,x),
\quad t\geq 0,\; x\in [0,1],\\
\ds u(0,x)&=&u_0(x) , \quad x\in[0,1],\\
\partial_x u(t,0) &=&\partial_x u(t,1)=0,\quad t> 0.
\end{array}
\right. \hskip-0.3cm
\end{equation}
Here $b,\sigma: \rr\mapsto \rr$ are the drift and diffusion coefficients
and $u_0:[0,1]\mapsto \rr$ is the initial condition. We write formally
$W(dt,dx)=\dot W(t,x)dtdx$, for $W(dt,dx)$ a white noise on 
$[0,\infty)\times [0,1]$ based on
$dtdx$, see Walsh \cite{w}.  
We will always assume in this paper that $b,\sigma$ are Lipschitz-continuous,
that is for some $C$,

\renewcommand{\theequation}{$\cH$}
\begin{equation} 
\hbox{for all $r,z \in \rr$, } \quad
|b(r)-b(z)|+|\sigma(r)-\sigma(z)| \leq C |r-z|.
\end{equation}
\renewcommand{\theequation}{\arabic{equation}}
\addtocounter{equation}{-1}

Our goals in this paper are the following:

$\bullet$ prove a strong existence and (partial)
uniqueness result when the initial condition 
$u_0$ only belongs to $L^1([0,1])$ and some stability results of the
solution with respect to such an initial condition;

$\bullet$ study the uniqueness of invariant measures and the asymptotic
confluence of solutions.

\vip

We will investigate these two points by using some {\it a priori} 
estimates on the difference between two solutions
$u,v$, obtained as a martingale dissipation of the 
$L^1([0,1])$-norm of $u(t)-v(t)$.

\vip

Let us mention that our results extend without difficulty
to the case of Dirichlet boundary conditions and to the case of 
the unbounded domain $\rr$ (with $u_0 \in L^1(\rr)$).

\vip

This equation has been much investigated, in particular 
since the work of Walsh \cite{w}. In \cite{w}, one can find definitions of 
weak solutions, existence and uniqueness results, as well as
proofs that solutions are H\"older-continuous, enjoy a Markov property, etc.
Let us mention for example
the works of Bally-Gyongy-Pardoux \cite{bgp} (existence of solutions
when the drift is only measurable), Gatarek-Goldys \cite{gg} (existence
of solutions in law), Donati-Pardoux (comparison results and reflection
problems), Bally-Pardoux (smoothness of the law of the solution),
Bally-Millet-Sanz \cite{bms} (support theorem), etc.
Sowers \cite{s}, Mueller \cite{m} and Cerrai \cite{c} have obtained some results
on the invariant distributions and convergence to equilibrium.

\subsection{Weak solutions}

We will consider two types of {\it weak} solutions, which we now precisely 
define, following the ideas of Walsh \cite{w}. When we refer to predictability,
this is with respect to the filtration $(\cF_t)_{t\geq 0}$ 
generated by $W$, that is $\cF_t=\sigma(W(A),A\in \cB([0,t]\times[0,1]))$.

\vip

We denote by $L^p([0,1])$ the set of all measurable functions 
$f:[0,1]\mapsto \rr$ such that $||f||_{L^p([0,1])}= 
(\int_0^1 |f(x)|^pdx)^{1/p}<\infty$.

\vip

Finally, we denote by $G_t(x,y)$ the Green kernel associated with the heat 
equation $\partial_t u = \partial_{xx}u$ on $\rr_+\times [0,1]$ 
with Neumann boundary conditions, whose explicit form can be found in 
Walsh \cite{w}.
Here we will only use that for some $C_T$, for all $x,y\in [0,1]$, 
all $t\in [0,T]$, see \cite{w},
\begin{equation}\label{ineqgt}
0\leq G_t(x,y)\leq \frac{C_T}{\sqrt t}e^{-|x-y|^2/4t}.
\end{equation}

\begin{defi} Assume $(\cH)$, and consider a $\rr$-valued
predictable process $u=(u(t,x))_{t\geq 0, x\in [0,1]}$.

(i) For $u_0\in L^1([0,1])$, $u$ 
is said to be a {\bf weak} solution to (\ref{she}) starting from
$u_0$ if a.s.,
\begin{equation}\label{bd1}
\hbox{for all } T>0, \quad
\sup_{[0,T]} ||u(t)||_{L^1([0,1])} 
+ \int_0^T ||\sigma(u(t))||_{L^2([0,1])}^2 dt<\infty
\end{equation}
and if for all $\varphi \in C^2_b([0,1])$ such
that $\varphi'(0)=\varphi'(1)=0$, for all $t\geq 0$,
a.s.,
\begin{eqnarray}\label{sheweak}
\int_0^1 u(t,x)\varphi(x) dx &=& \int_0^1 u_0(x)\varphi(x) dx
+ \int_0^t \int_0^1 \sigma(u(s,x))\varphi(x) W(ds,dx) \\
&&+ \int_0^t \int_0^1 [u(s,x)\varphi''(x) + b(u(s,x))\varphi(x)]dxds.
\nonumber
\end{eqnarray}
(ii) For $u_0$ bounded-measurable, $u$
is said to be a {\bf mild} solution to (\ref{she}) starting from
$u_0$ if a.s.,
\begin{equation}\label{bd2}
\hbox{for all } T>0, \quad
\sup_{[0,T]\times [0,1]} |u(t,x)|<\infty
\end{equation}
and if for all $t\geq 0$, all $x\in [0,1]$, a.s.,
\begin{equation}\label{shemild}
u(t,x)= \int_0^1 G_t(x,y)u_0(y)dy + \intot \int_0^1 G_{t-s}(x,y)
[\sigma(u(s,y)) W(ds,dy)+ b(u(s,y)) dyds]. 
\end{equation}
\end{defi}

Let us make a few comments.
Recall that for $(H(s,y))_{s\geq 0, y\in [0,1]}$ a $\rr$-valued 
predictable process,
the stochastic integral $\int_0^t \int_0^1 H(s,y)W(ds,dy)$ is well-defined
if and only if $\int_0^t \int_0^1 H^2(s,y) dyds <\infty$ a.s.

\vip

$\bullet$ Thus (\ref{bd1}) implies that all the terms in 
(\ref{sheweak})
are well-defined. Clearly, condition (\ref{bd1}) is not far from minimal.

\vip

$\bullet$ Next, (\ref{bd2}) and (\ref{ineqgt}) imply that all the terms 
in (\ref{shemild})
are well-defined, but here (\ref{bd2}) is clearly far from optimal.

\vip

When $u_0$ only belongs to $L^1([0,1])$, we will only be able to prove that 
(\ref{bd1}) holds.

\vip

Let us finally recall that Walsh \cite{w} proved, under $(\cH)$, that for any
bounded-measurable initial condition $u_0$, there exists a unique
mild solution $u$ to (\ref{she}), 
which is also a weak solution and which furthermore
satisfies, for all $p\geq 1$, all $T>0$,
$\E[\sup_{[0,T]\times[0,1]} |u(t,x)|^p ]<\infty$.

\subsection{Existence and stability in $L^1([0,1])$}

Our first goal is to extend the existence theory to more
general initial conditions.

\begin{theo}\label{exstab}
Assume $(\cH)$.

(i) For $u_0 \in L^1([0,1])$, there exists
a weak solution $u$ 
to (\ref{she}) starting from $u_0$.

(ii) This solution is unique in the following sense: for any sequence
of bounded-measurable functions $u_0^n:[0,1]\mapsto \rr$ such that
$\lim_n ||u_0^n-u_0||_{L^1([0,1])}=0$, the sequence 
$\sup_{[0,T]} ||u^n(t)-u(t)||_{L^1([0,1])}$ tends to $0$ in probability for
any $T$. Here $u^n$ is the unique mild solution to (\ref{she})
starting from $u_0^n$.

(iii) For $u_0,v_0 \in L^1([0,1])$, consider the two weak solutions 
$u$ and $v$ to (\ref{she}) starting from $u_0$ and $v_0$ built in (i). 
For all $\gamma \in (0,1)$, all $T\geq 0$, we have
$$
\E\left[\sup_{[0,T]}||u(t)-v(t)||_{L^1([0,1])}^\gamma + 
\left(\int_0^T ||\sigma(u(t))-\sigma(v(t))||_{L^2([0,1])}^2dt \right)^{\gamma/2}  
\right] \leq C_{b,\gamma,T} ||u_0-v_0||_{L^1([0,1])}^\gamma,
$$
where $C_{b,\gamma,T}$ depends only on $b,\gamma,T$.

(iv) Assume now that $b$ is non-increasing.  
For $u_0,v_0 \in L^1([0,1])$, let $u,v$ be the two weak solutions 
to (\ref{she}) starting from $u_0$ and $v_0$ built in (i).
For all $\gamma \in (0,1)$, we have
\begin{eqnarray*}
\E\Big[\sup_{[0,\infty)}||u(t)-v(t)||_{L^1([0,1])}^\gamma +
\left(\int_0^\infty ||b(u(t))-b(v(t)) ||_{L^1([0,1])} dt \right)^\gamma &&\\ 
+\left(\int_0^\infty 
||\sigma(u(t))-\sigma(v(t))||_{L^2([0,1])}^2dt \right)^{\gamma/2}  
\Big] &\leq& C_{\gamma} ||u_0-v_0||_{L^1([0,1])}^\gamma,
\end{eqnarray*}
where $C_{\gamma}$ depends only on $\gamma$.
\end{theo}

Observe that this result contains a regularization property.
For example if $\sigma(z)=z$,
even if $u_0$ does not belong to $L^2([0,1])$, the weak solution 
satisfies (\ref{bd1}) and in particular $\sigma(u(t))=u(t)\in L^2([0,1])$ 
for a.e. $t>0$. For the same reasons, the stability result (iii)
provides a better estimate for a.e. $t>0$ than for $t=0$.

\vip

To our knowledge, Theorem \ref{exstab} is the first result concerning
$L^1([0,1])$ initial conditions. Many works concern bounded-measurable 
(or continuous)
initial conditions, see Walsh \cite{w}, Bally-Gyongy-Pardoux \cite{bgp},
Cerrai \cite{c}. Another abundant literature deals with the Hilbert
case (initial conditions in $L^2([0,1])$), see Pardoux \cite{p}, 
Da Prato-Zabczyk
\cite{dz}, Gatarek-Goldys \cite{gg}.

\vip

The present well-posedness result is quite satisfying, since the requirement
that $u_0 \in L^1([0,1])$ is very weak and seems necessary for 
(\ref{sheweak}) to make sense.

\subsection{Large time behavior}

We now wish to study the uniqueness of invariant measures. 

\begin{defi}\label{dinv}
A probability measure $Q$ on $L^1([0,1])$ is said to be an invariant
distribution for (\ref{she}) if, for $u_0$ a $L^1([0,1])$-valued
random variable with law $Q$ independent of $W$, for $u$
the weak solution to (\ref{she}) starting from $u_0$ built
in Theorem \ref{exstab}, $\cL(u(t))=Q$ for all $t\geq 0$.
\end{defi}

We have the following result.

\begin{theo}\label{uniinv}
Assume $(\cH)$, that $b$ is non-increasing and that 
$(\sigma,b):\rr\mapsto\rr^2$ is injective. Then
(\ref{she}) admits at most one invariant distribution.
\end{theo}

To prove the asymptotic confluence of solutions, 
we need to strengthen the injectivity assumption.

\renewcommand{\theequation}{$\cI$}
\begin{equation}
\left\{
\begin{array}{l}
\hbox{There is a strictly increasing convex function } 
\rho:\rr_+\mapsto \rr_+ \hbox{ with } \rho(0)=0 \hbox{ such that }\\
\hbox{for all $r,z\in \rr$, } \quad
|b(r)-b(z)|+|\sigma(r)-\sigma(z)|^2 \geq \rho(|r-z|).
\end{array}
\right.
\end{equation}
\renewcommand{\theequation}{\arabic{equation}}
\addtocounter{equation}{-1}

\begin{theo}\label{confl}
Assume $(\cH)$, that $b$ is non-increasing and $(\cI)$.

(i) The following asymptotic 
confluence property holds: for $u_0,v_0\in L^1([0,1])$, 
for $u,v$
the weak solutions to (\ref{she}) starting from $u_0$ and $v_0$
built in Theorem \ref{exstab},
$$
\hbox{a.s., } \quad \lim_{t\to \infty} ||u(t)-v(t)||_{L^1([0,1])} =0.
$$
(ii) Assume additionally that (\ref{she}) admits an invariant 
distribution $Q$. Then for $u_0\in L^1([0,1])$, for
$u$ the corresponding weak solution
to (\ref{she}), $u(t)$ goes in law to $Q$ as $t\to\infty$.
\end{theo}

Clearly, $(\cI)$ holds if $b$ is $C^1$ with $b'\leq -\e<0$ 
(choose $\rho(z)=\e z$)
or if $\sigma$ is $C^1$ with $|\sigma'|\geq \e>0$ (choose $\rho(z)=(\e z)^2$).
One may also combine conditions on $b$ and $\sigma$.

But $(\cI)$ also holds if $b$ is $C^1$ and if $b'\leq 0$ vanishes reasonably.
For example if  $b(z)=-\sg(z)\min(|z|,|z|^p)$ for some $p\geq 1$,
choose $\rho= \e \rho_p$ with $\e$ small enough and  
$\rho_p(z)= z^p$ for $z\in [0,1]$ 
and $\rho_p(z)= p z -p+ 1$ for $z\geq 1$. If 
$b(z)=-z-\sin z$, choose $\rho=\e \rho_3$ with $\e$ small enough.

One may also consider the case where 
$\sigma$ is monotonous with $\sigma'$ vanishing 
reasonably.
\vip

Let us now compare Theorems \ref{uniinv} and \ref{confl} with known results.
The works cited below sometimes
concern different boundary conditions, but we believe this is not
important.

\vip

$\bullet$ Sowers \cite{s} has proved the existence of an invariant
distribution supported by $C([0,1])$, assuming $(\cH)$, 
that $\sigma$ is bounded
and that $b$ is of the form $b(z)=- \alpha z + f(z)$, for some bounded $f$
and some $\alpha>0$.
He obtained uniqueness of this invariant distribution
when $\sigma$ is sufficiently small and bounded from below.

\vip

$\bullet$ Mueller \cite{m} has obtained some surprising coupling results,
implying in particular the uniqueness of
an invariant distribution as well as a the trend to equilibrium.
He assumes 
$(\cH)$, that $\sigma$ is bounded from above and from below and that
$b$ is non-increasing, with $|b(z)-b(r)|\geq \alpha |z-r|$ 
for some $\alpha>0$.

\vip

$\bullet$ Cerrai \cite{c} assumed that $\sigma$ is strictly
monotonous (it may vanish, but only at one point). 

(i) She obtained an asymptotic 
confluence result which we do not recall here and concerns, 
roughly, the case
$b(z)\simeq -\sg(z) |z|^{m}$ as $z\to \pm \infty$, for some $m>1$.

(ii) Assuming $(\cH)$, she proved uniqueness of the invariant
distribution as well as an asymptotic confluence property, under the conditions
that for all $r\leq z$, $b(z)-b(r)\leq \lambda (z-r)$,
and $|\sigma(z)-\sigma(r)| \geq \mu |z-r|$, for some $\mu>0$
and some $\lambda < \mu^2/2$
(if $b$ is non-increasing, choose $\lambda=0$).

\vip

Thus the main advantages of the present paper are that  
the uniqueness of the invariant measure requires very few
conditions, and we allow $\sigma$
to vanish (it may be compactly supported).

\vip



{\it Example 1.} Assume $(\cH)$ and that $b$ strictly decreasing.
Then there exists at most one invariant distribution.
If $b(z)=-z$ or $b(z)=-z-\sin z$ or $b(z)=-\sg(z)
\min(|z|,|z|^p)$ for some $p>1$,
then we have  asymptotic confluence of solutions.
Here to apply \cite{s,m} one needs to assume additionally that
$\sigma$ is bounded from above and from below, while to apply
\cite{c}, one has to suppose that $\sigma$ is strictly monotonous.

\vip

{\it Example 2.} Assume $(\cH)$, that $b$ is non-increasing and that
$\sigma$ is strictly monotonous.
Then there exists at most one invariant distribution.

If furthermore $\sigma$ is $C^1$ with $0< c< \sigma'<C$, 
then we get  asymptotic 
confluence of solutions using \cite{c} or Theorem \ref{confl}
(here \cite{s,m} cannot apply, since $\sigma$ vanishes).
But now if $\sigma'\geq 0$ reasonably vanishes 
then  Theorem \ref{confl} applies, which is not the case of \cite{c}:
take e.g. $\sigma(z)= \sg(z) \min(|z|,|z|^p)$ for some $p>1$,
or $\sigma(z)=z+\sin z$.

\vip

{\it Example 3.} Consider the compactly supported coefficient
$\sigma(z)=(1-z^2)\indiq_{\{|z|\leq 1\}}$.
Assume that $b$ is $C^1$, non-increasing, with $b'(z) \leq -\e<0$ for
$z\in (-\infty,-1)\cup\{0\}\cup(1,+\infty)$. 
Then Theorems \ref{uniinv} and \ref{confl} apply, while \cite{s,m,c} do not.

Observe here that if $b(z_0)=0$ for some $z_0\notin (-1,1)$, then
$u(t)\equiv z_0$ is the (unique) stationary solution.

If now $b(-1)>0$ and $b(1)<0$, 
then the invariant measure $Q$
(that exists due to Sowers \cite{s}) is unique and one may show, 
using the comparison
Theorem of Donati-Pardoux \cite{dp}, that
$Q$ is supported by $[-1,1]$-valued continuous functions on $[0,1]$. 

\vip

However, there are some cases where \cite{s,c} provide some better
results than ours.

\vip

{\it Example 4.} If $\sigma(z)=\mu z$ and $b(z)=\lambda z$, then
$u(t)\equiv 0$ is an obvious stationary solution.
Theorems \ref{uniinv} and \ref{confl} apply if 
$\lambda \leq 0$ and $|\lambda|+|\mu|>0$. Cerrai \cite{c} was
able to treat the case $\lambda>0$ provided $\mu^2/2>\lambda$.

\vip

{\it Example 5.} If $\sigma$ is small enough and bounded from below
and if $b(z)=-\alpha z + h(z)$, with $\alpha>0$ and $h$ bounded, 
then Sowers \cite{s} obtains the uniqueness of the invariant distribution
even if $b$ is not non-increasing.


\subsection{Plan of the paper}
In the next section, we prove some inequalities concerning
the $L^1([0,1])$-norm of the difference between any pair of {\it mild} solutions
to (\ref{she}). Section \ref{exist} is dedicated to the proof
of our existence result Theorem \ref{exstab}.
Theorems \ref{uniinv} and \ref{confl} are checked in Section \ref{ergo}.
We briefly discuss the multi-dimensional equation in Section \ref{dim2}
and conclude the paper with an appendix containing technical results.

\section{On the $L^1([0,1])$-norm of the difference
between two mild solutions}\label{basic}

All our study is based on the following result. We set
$\sg(z)=1$ for $z\geq 0$ and $\sg(z)=-1$ for $z< 0$.

\begin{pro}\label{eql1}
Assume $(\cH)$. 
For two bounded-measurable initial conditions 
$u_0,v_0$, let $u,v$ be the corresponding mild
solutions to (\ref{she}). Then, enlarging the probability space if necessary,
there is a 
Brownian motion $(B_t)_{t\geq 0}$ such that a.s., for all
$t\geq 0$, 
\begin{eqnarray}\label{ineqfc}
||u(t)-v(t)||_{L^1([0,1])}&\leq&||u_0-v_0||_{L^1([0,1])}
+ \intot || \sigma(u(s))-\sigma(v(s)) ||_{L^2([0,1])} dB_s \\
&&+\intot  \int_0^1  \sg(u(s,x)-v(s,x))(b(u(s,x))-b(v(s,x)))dxds.\nonumber
\end{eqnarray}
\end{pro}

\begin{proof} We divide the proof into several steps, following closely
the ideas of Donati-Pardoux \cite[Theorem 2.1]{dp}, to which we refer
for technical details.

\vip

{\it Step 1.} Consider an orthonormal basis $(e_k)_{k\geq 1}$ of
$L^2([0,1])$. For $k\geq 1$, we 
set $B^k_t=\int_0^t \int_0^1 e_k(x) W(ds,dx)$. Then
$(B^k)_{k\geq 1}$ is a family of independent Brownian motions.
For $n\geq 1$, consider the unique adapted solution 
$u^n\in L^2(\Omega\times[0,T],V)$, where $V=\{f\in H^1([0,1]), f'(0)=f'(1)=0\}$,
to
$$
u^n(t,x)=u_0(x)+\intot \left[\partial_{xx}u^n(s,x)ds + b(u^n(s,x))\right]ds 
+ \sum_{k=1}^n \int_0^t \sigma(u^n(s,x))e_k(x)dB^k_s.
$$
We refer to Pardoux \cite{p} for existence, uniqueness and properties
of this solution. We also consider the solution $v^n$ to the same equation
starting from $v_0$. Then, as shown in \cite{dp},
\begin{equation}\label{cacv}
\lim_n \sup_{[0,T]\times[0,1]}\E[|u^n(t,x)-u(t,x)|^2 +|v^n(t,x)-v(t,x)|^2 ]=0.
\end{equation}

{\it Step 2.} For $\e>0$, we introduce a nonnegative
$C^2$ function $\phi_\e$ such that $\phi_\e(z)= |z|$ for $|z|\geq \e$,
with $|\phi_\e'(z)| \leq 1$ and $0\leq \phi_\e''(z)
\leq 2\e^{-1}\indiq_{|z|<\e}$.  
When applying the It\^o formula (see \cite{dp} for details), we get
\begin{eqnarray}\label{pbfc}
&&\hskip-1.5cm\int_0^1 \phi_\e(u^n(t,x)-v^n(t,x))dx
= \int_0^1 \phi_\e(u_0(x)-v_0(x))dx \\
&&+ \intot \int_0^1 \phi_\e'(u^n(s,x)-v^n(s,x))
\partial_{xx}[u^n(s,x)-v^n(s,x)]dxds \ala
&&+ \intot \int_0^1 \phi_\e'(u^n(s,x)-v^n(s,x))[b(u^n(s,x))-b(v^n(s,x))]dxds\ala
&&+ \sum_{k=1}^n \int_0^t \int_0^1 \phi_\e'(u^n(s,x)-v^n(s,x)) 
[\sigma(u^n(s,x))-\sigma(v^n(s,x)) ] e_k(x)dx dB^k_s\ala
&&+ \frac{1}{2}  \sum_{k=1}^n\intot \int_0^1 \phi_\e''(u^n(s,x)-v^n(s,x)) 
[\sigma(u^n(s,x))-\sigma(v^n(s,x)) ]^2 e^2_k(x)dx ds\ala
&&=:I_\e^1+I_\e^2(t)+I_\e^3(t)+I_\e^4(t)+I_\e^5(t).\nonumber
\end{eqnarray}
Since $|z|\leq \phi_\e(z)\leq |z|+\e$ for all $z$, we easily get, a.s.,
$$
\lim_{\e\to 0}\int_0^1 \phi_\e(u^n(t,x)-v^n(t,x))dx =||u^n(t)-v^n(t)||_{L^1([0,1])}
\quad \hbox{ and } \quad
\lim_{\e\to 0} I^1_\e =||u_0-v_0||_{L^1([0,1])}.
$$
An integration by parts, using that $\partial_x[u^n(t,0)-v^n(t,0)]=
\partial_x[u^n(t,1)-v^n(t,1)]=0$ shows that
$$
I^2_\e(t)=-\intot \int_0^1 \phi_\e''(u^n(s,x)-v^n(s,x))
[\partial_{x}(u^n(s,x)-v^n(s,x))]^2\leq 0.
$$
Since $\phi_\e''(z-r)(\sigma(z)-\sigma(r))^2\leq C 
\e^{-1} \indiq_{|z-r|\leq \e} |z-r|^2 \leq C \e$ by $(\cH)$, 
we have $I^5_\e(t) \leq Cnt\e $, whence   
$$
\lim_{\e\to 0} I_\e^5(t)=0 \hbox{ a.s.}
$$
Using that $|\phi_\e'(z)- \sg(z)|\leq \indiq_{\{|z|\leq \e\}}$ and $(\cH)$,
one obtains a.s.
\begin{eqnarray*}
&&\lim_{\e\to 0}\left| I^3_\e(t) - \int_0^t\int_0^1  \sg(u^n(s,x)-v^n(s,x))
(b(u^n(s,x))-b(v^n(s,x))) dxds \right| \\
&\leq & \lim_{\e\to 0} \int_0^t \int_0^1 \indiq_{|u^n(s,x)-v^n(s,x)|\leq \e}
|b(u^n(s,x))-b(v^n(s,x))| dxds  \leq \lim_{\e\to 0} Ct \e =0.
\end{eqnarray*}
Similarly,
\begin{eqnarray*}
&&\lim_{\e\to 0} \E\left[\left(I^4_\e(t)-
\sum_{k=1}^n \int_0^t \int_0^1 \sg(u^n(s,x)-v^n(s,x))
[\sigma(u^n(s,x))-\sigma(v^n(s,x)) ] e_k(x)dx dB^k_s \right)^2 \right]=0.
\end{eqnarray*}
Thus we can pass to the limit as $\e\to 0$ in (\ref{pbfc}) and get, a.s.,
\begin{eqnarray}\label{pbfcf}
||u^n(t)-v^n(t)||_{L^1([0,1])}&\leq&||u_0-v_0||_{L^1([0,1])} \nonumber \\
&&+\intot  \int_0^1  \sg(u^n(s,x)-v^n(s,x))[b(u^n(s,x))-b(v^n(s,x))]dxds
\nonumber\\
&&+ \sum_{k=1}^n \intot \int_0^1
\sg(u^n(s,x)-v^n(s,x)) [\sigma(u^n(s,x))-\sigma(v^n(s,x))] 
e_k(x)dx dB^k_s .
\end{eqnarray}

{\it Step 3.} Using $(\cH)$, there holds, for all $r_1,z_1,r_2,z_2$ in
$\rr$,
\begin{eqnarray}\label{sglip}
\big|\sg(r_1-z_1)[\sigma(r_1)-\sigma(z_1)]
-\sg(r_2-z_2)[\sigma(r_2)-\sigma(z_2)]  \big| 
\leq C(|r_1-r_2|+|z_1-z_2|),\\
\big|\sg(r_1-z_1)[b(r_1)-b(z_1)]
-\sg(r_2-z_2)[b(r_2)-b(z_2)]  \big| 
\leq C(|r_1-r_2|+|z_1-z_2|).
\end{eqnarray}
Indeed, it suffices, by symmetry, to check that 
$\big|\sg(r_1-z_1)[\sigma(r_1)-\sigma(z_1)]
-\sg(r_2-z_1)[\sigma(r_2)-\sigma(z_1)]  \big| 
\leq C|r_1-r_2|$. If $\sg(r_1-z_1)=\sg(r_2-z_2)$, this is obvious.
If now $r_1\leq z_1 \leq r_2$ (or $r_1\geq z_1 \geq r_2$) 
we get the upper-bound
$|\sigma(r_1)+\sigma(r_2)-2\sigma(z_1)| \leq C(|r_1-z_1|+|r_2-z_1|)
= C |r_1-r_2|$.

Using (\ref{cacv}), it is thus routine
to make $n$ tend to infinity in (\ref{pbfcf}) and to obtain, a.s.,
\begin{eqnarray}\label{coucou}
||u(t)-v(t)||_{L^1([0,1])}&\leq&||u_0-v_0||_{L^1([0,1])} 
+\intot  \int_0^1  \sg(u(s,x)-v(s,x))[b(u(s,x))-b(v(s,x))]dxds \nonumber \\
&&+ \sum_{k=1}^\infty \intot \int_0^1
\sg(u(s,x)-v(s,x)) [\sigma(u(s,x))-\sigma(v(s,x))] 
e_k(x)dx dB^k_s. 
\end{eqnarray}
For the last term, we used that, by the Plancherel identity, 
setting for simplicity 
\begin{eqnarray*}
\alpha_n(s,x)&=&\sg(u^n(s,x)-v^n(s,x)) [\sigma(u^n(s,x))-\sigma(v^n(s,x))],\\
\alpha(s,x)&=&\sg(u(s,x)-v(s,x)) [\sigma(u(s,x))-\sigma(v(s,x))],
\end{eqnarray*}
there holds
\begin{eqnarray*}
&&\E\Big[\Big( \sum_{k=1}^n \intot \int_0^1
\alpha_n(s,x) e_k(x) dB^k_s - \sum_{k=1}^\infty \intot \int_0^1
\alpha(s,x) e_k(x) dB^k_s\Big)^2 \Big]\ala
&\leq& \intot \E\Big[ \sum_{k\geq 1} \Big(\int_0^1\Big\{ 
\alpha_n(s,x)-\alpha(s,x)\Big\}e_k(x)dx \Big)^2 \Big]ds
+ \sum_{k\geq n+1} \intot \E\Big[ \Big(\int_0^1 \alpha(s,x)e_k(x)dx \Big)^2 
\Big]ds\\
&\leq& \intot \E\Big[ || \alpha_n(s)-\alpha(s)||_{L^2([0,1])}^2\Big]ds 
+ \sum_{k\geq n+1} \intot \E\Big[ \Big(\int_0^1
\alpha(s,x) e_k(x)dx  \Big)^2 \Big]ds=:I_n(t)+J_n(t).
\end{eqnarray*}
Using (\ref{sglip}) and then (\ref{cacv}), $I_n(t)\leq
C \int_0^t \int_0^1 \E[|u^n(s,x)-u(s,x)|^2+|v^n(s,x)-v(s,x)|^2 ]dxds$ 
tends to $0$ as $n\to \infty$. Finally, $J_n(t)$ tends to $0$
because $\sum_{k\geq 1}  \int_0^t \E [(\int_0^1
\alpha(s,x) e_k(x)dx )^2 ] ds  = \int_0^t \E[||\alpha(s)||_{L^2([0,1])}^2] ds
\leq C\int_0^t \int_0^1 \E(|u(s,x)-v(s,x)|^2) dxds <\infty$.

\vip

{\it Step 4.}
A standard representation argument (see e.g. Revuz-Yor 
\cite[Proposition 3.8 and Theorem 3.9 p 202-203]{ry}) 
concludes the proof, because the last term on the RHS of
(\ref{coucou}) is a continuous
local martingale with bracket 
\begin{equation*}
\intot \sum_{k=1}^\infty \left(\int_0^1
\sg(u(s,x)-v(s,x))[\sigma(u(s,x))-\sigma(v(s,x))]e_k(x)dx \right)^2 ds 
=\intot ||\sigma(u(s))-\sigma(v(s))||_{L^2([0,1])}^2 ds.
\end{equation*}
We used here again that $(e_k)_{k\geq 1}$ is an orthonormal 
basis of $L^2([0,1])$.
\end{proof}


\begin{cor}\label{ineql1g}
Adopt the notation and assumptions of Proposition \ref{eql1}.
For all $\gamma \in (0,1)$, all $T\geq 0$, 
$$
\E\left[\sup_{[0,T]}||u(t)-v(t)||_{L^1([0,1])}^\gamma + 
\left(\int_0^T ||\sigma(u(t))-\sigma(v(t))||_{L^2([0,1])}^2dt \right)^{\gamma/2}  
\right] \leq C_{b,\gamma,T} ||u_0-v_0||_{L^1([0,1])}^\gamma,
$$
where $C_{b,\gamma,T}$ depends only on $b,\gamma,T$.
\end{cor}

\begin{proof}
Let $C$ be the Lipschitz constant of $b$. Denote by $L_t$ the RHS of 
(\ref{ineqfc}). The It\^o formula yields 
\begin{eqnarray*}
||u(t)-v(t)||_{L^1([0,1])}e^{-Ct}&\leq& L_te^{-Ct} \\
&=& ||u_0-v_0||_{L^1([0,1])}- C \intot  
e^{-Cs }L_s ds\\
&&+ \intot || \sigma(u(s))-\sigma(v(s)) ||_{L^2([0,1])}e^{-Cs} dB_s \\
&&+\intot  \int_0^1  e^{-Cs }\sg(u(s,x)-v(s,x))(b(u(s,x))-b(v(s,x)))dxds.
\end{eqnarray*}
But $\int_0^1 \sg(u(s,x)-v(s,x))(b(u(s,x))-b(v(s,x)))dx \leq
C ||u(s)-v(s)||_{L^1([0,1])} \leq C L_s$.
Hence
\begin{eqnarray*}
||u(t)-v(t)||_{L^1([0,1])}e^{-Ct}&\leq& ||u_0-v_0||_{L^1([0,1])}
+\intot || \sigma(u(s))-\sigma(v(s)) ||_{L^2([0,1])}e^{-Cs} dB_s=:M_t.
\end{eqnarray*}
Hence $M_t$ is a nonnegative local martingale with bracket 
$\lb M \rb_t= \int_0^t || \sigma(u(s))-\sigma(v(s)) ||_{L^2([0,1])}^2
e^{-2Cs} ds$. Applying Lemma \ref{martingales}, we immediately get,
for $\gamma \in (0,1)$,
\begin{eqnarray*}
&&\E\left[\sup_{[0,\infty)}||u(t)-v(t)||_{L^1([0,1])}^\gamma e^{-C\gamma t}  
+ \left( \int_0^\infty || \sigma(u(s))-\sigma(v(s)) ||_{L^2([0,1])}^2
e^{-2Cs} ds \right)^{\gamma/2} \right] \\
&\leq& C_\gamma ||u_0-v_0||_{L^1([0,1])}^\gamma.
\end{eqnarray*}
The result easily follows.
\end{proof}

Finally, one can say a little more when $b$ is non-increasing.

\begin{cor}\label{ineql1d}
Adopt the notation and assumptions of Proposition \ref{eql1} and assume
that $b$ is non-increasing. Then for all $\gamma \in (0,1)$,
\begin{eqnarray*}
\E\Big[\sup_{[0,\infty)}||u(t)-v(t)||_{L^1([0,1])}^\gamma + 
\left(\int_0^\infty ||b(u(t))-b(v(t))||_{L^1([0,1])}dt \right)^{\gamma} &&\\ 
+\left(\int_0^\infty 
||\sigma(u(t))-\sigma(v(t))||_{L^2([0,1])}^2dt \right)^{\gamma/2}  
\Big] &\leq& C_{\gamma}||u_0-v_0||_{L^1([0,1])}^\gamma,\nonumber
\end{eqnarray*}
where $C_{\gamma}$ depends only on $\gamma$.
\end{cor}

\begin{proof}
Since $b$ is non-increasing, Proposition \ref{eql1} yields
\begin{eqnarray*}
&&||u(t)-v(t)||_{L^1([0,1])}+ \intot ||b(u(s))-b(v(s))||_{L^1([0,1])}ds \\
&\leq& ||u_0-v_0||_{L^1([0,1])}
+ \intot || \sigma(u(s))-\sigma(v(s)) ||_{L^2([0,1])} dB_s=:M_t,
\end{eqnarray*}
which is thus a nonnegative martingale with bracket
$\lb M \rb_t= \int_0^t || \sigma(u(s))-\sigma(v(s)) ||_{L^2([0,1])}^2 ds$.
Lemma \ref{martingales} allows us to conclude.
\end{proof}

\section{Existence theory in $L^1([0,1])$}\label{exist}

The goal of this section is to give the

\vip

\begin{preuve} {\it of Theorem \ref{exstab}.} We start with point (i).
Let thus $u_0\in L^1([0,1])$ and consider a sequence of bounded-measurable
initial conditions $(u_0^n)_{n\geq 1}$ 
such that $||u_0^n-u_0||_{L^1([0,1])}\leq 2^{-n}$. For each $n\geq 1$, 
denote by $u^n$ the mild solution
to (\ref{she}) starting from $u_0^n$. Using Corollary \ref{ineql1g}
(with $\gamma=1/2$), we deduce that a.s.,
$$
\sum_{n\geq 1} \left[\sup_{[0,T]}||u^{n+1}(t)-u^n(t)||_{L^1([0,1])}^{1/2}
+ \left(\int_0^T ||\sigma(u^{n+1}(t))-\sigma(u^n(t))||_{L^2([0,1])}^2   dt
\right)^{1/4} \right]< \infty,
$$
which implies that
$$
\sum_{n\geq 1} \left[\sup_{[0,T]}||u^{n+1}(t)-u^n(t)||_{L^1([0,1])} 
+ ||\sigma(u^{n+1})-\sigma(u^n)||_{L^2([0,T]\times[0,1])} \right]< \infty.
$$
Using some completeness arguments, we deduce that there are some
(predictable) processes $u$ and $S$ such that a.s., for all $T>0$,
$\sup_{[0,T]}||u(t)||_{L^1([0,1])}+\int_0^T ||S(t)||_{L^2([0,1])}^2dt <\infty$ and
$$
\lim_n \sup_{[0,T]}||u(t)-u^n(t)||_{L^1([0,1])} = 0, \quad
\lim_n ||S-\sigma(u^n)||_{L^2([0,T]\times[0,1])} = 0.
$$
Since $\sigma$ is Lipschitz-continuous, we deduce from 
the first equality that $\lim_n ||\sigma(u)-\sigma(u^n)||_{L^1([0,T]\times[0,1])}
= 0$, 
while from the second one, 
$\lim_n ||S-\sigma(u^n)||_{L^1([0,T]\times[0,1])}= 0$. Consequently, $S=\sigma(u)$
a.e. and we finally conclude that a.s.,
\begin{equation}\label{ettac}
\hbox{ for all } T>0,\quad
\lim_n \left(\sup_{[0,T]}||u(t)-u^n(t)||_{L^1([0,1])} +\int_0^T ||\sigma(u(t)) 
-\sigma(u^n(t))||_{L^2([0,1])}^2 dt\right)= 0.
\end{equation}
It remains to prove that $u$ is a weak solution to (\ref{she}).
We have already seen that $u$ satisfies (\ref{bd1}). Next,
for $\varphi\in C^2_b([0,1])$ with $\varphi'(0)=\varphi'(1)=0$, for
$t\geq 0$, we know that a.s., $A^{n,\varphi}_t=B^{n,\varphi}_t$ for all
$n\geq 1$, where
\begin{eqnarray*}
A^{n,\varphi}_t&:=&\int_0^1 \varphi(x)u^n(t,x)dx 
- \int_0^1 \varphi(x)u^n_0(x)dx  
- \intot \int_0^1 [u^n(s,x)\varphi''(x) + b(u^n(s,x))\varphi(x)] dxds
\\
B^{n,\varphi}_t&:=& \intot \int_0^1 \sigma(u^n(s,x))\varphi(x) W(ds,dx).
\end{eqnarray*}
It directly follows from (\ref{ettac}) and $(\cH)$ that a.s.,
$$
\lim_{n\to \infty} A^{n,\varphi}_t=\int_0^1 \varphi(x)u(t,x)dx 
- \int_0^1 \varphi(x)u_0(x)dx  
- \intot \int_0^1 [u(s,x)\varphi''(x) + b(u(s,x))\varphi(x)] dxds.
$$
We deduce that $B^{\varphi}_t:=\lim_n B^{n,\varphi}_t$ exists a.s.,
and it only remains to check that $B^\varphi_t=C^\varphi_t$ a.s., where
$C^\varphi_t:=\int_0^t \int_0^1 
\sigma(u(s,x))\varphi(x) W(ds,dx)$.
To this end, consider, for $M>0$,  the stopping time
$$
\tau_M = \inf\left\{ r \geq 0, \;  \int_0^r ||\sigma(u(s))||_{L^2([0,1])}^2ds   
+ \sup_n \int_0^r ||\sigma(u^n(s))||_{L^2([0,1])}^2ds  \geq M\right\}.
$$
Using (\ref{ettac}) and the dominated convergence Theorem, 
we see that for each $M>0$,
$$
\lim_n \E[|B^{n,\varphi}_{t\land \tau_M}- C^{\varphi}_{t\land \tau_M} |^2]
= \lim_n \E\left[\int_0^{t\land \tau_M} 
||(\sigma(u(s))-\sigma(u^n(s))\varphi ||_{L^2([0,1])}^2 ds \right] = 0.
$$
But we also deduce from (\ref{ettac}) that  a.s.,
$\sup_n \int_0^T ||\sigma(u^n(s))||_{L^2([0,1])}^2ds <\infty$ for all $T>0$,
whence $\lim_{M \to \infty} \tau_M = \infty$ a.s.
We easily conclude that $B^{n,\varphi}_t$ tends to $C^\varphi_t$ in probability,
whence $B^\varphi_t=C^\varphi_t$ a.s.

\vip

Point (ii) is easily checked: let $(\tu_0^n)_{n\geq 1}$ 
be another sequence of bounded-measurable 
initial conditions converging to $u_0$ and let
$(\tu^n)_{n\geq 1}$ be the corresponding sequence of mild solutions 
to (\ref{she}). Then necessarily, $||u_0^n - \tu_0^n||_{L^1([0,1])}$ 
tends to $0$,
whence, by Corollary \ref{ineql1g}, $\sup_{[0,T]}||u^n(t)-\tu^n(t)||_{L^1([0,1])}$
tends also to $0$, in probability. Using (\ref{ettac}), 
we conclude that $\sup_{[0,T]}||u(t)-\tilde u^n(t)||_{L^1([0,1])}$ tends to $0$
in probability.

\vip

We now prove point (iii). For $u_0$ and $v_0$ in $L^1([0,1])$, 
we consider $u_0^n$ and $v_0^n$ bounded-measurable with 
$||u_0^n-u_0||_{L^1([0,1])}
+||v_0^n-v_0||_{L^1([0,1])} \leq 2^{-n}$. We denote by $u,v,u^n,v^n$ the 
corresponding
weak solutions to (\ref{she}). In the proof of (i), we have seen that
a.s., $\lim_n \sup_{[0,T]} [||u^n(t)-u(t)||_{L^1([0,1])}
+||v^n(t)-v(t)||_{L^1([0,1])} ]=0$ and 
$\lim_n \int_0^T [||\sigma(u^n(t))-\sigma(u(t))||_{L^2([0,1])}^2
+||\sigma(v^n(t))-\sigma(v(t))||_{L^2([0,1])}^2 ]dt =0$. 
Using the Fatou Lemma and Corollary \ref{ineql1g}, we thus get
\begin{eqnarray*}
&&\E\left[\sup_{[0,T]}||u(t)-v(t)||_{L^1([0,1])}^\gamma 
+ \left(\int_0^T ||\sigma(u(t))-\sigma(v(t))||_{L^2([0,1])}^2 dt
\right)^{\gamma/2} \right] \ala
&\leq& \liminf_n \E\left[\sup_{[0,T]}||u^{n}(t)-v^n(t)||_{L^1([0,1])}^\gamma 
+ \left(\int_0^T ||\sigma(u^{n}(t))-\sigma(v^n(t))||_{L^2([0,1])}^2 dt   
\right)^{\gamma/2}  \right] \\
&\leq& \liminf_n C_{\gamma,T} ||u_0^n-v_0^n||_{L^1([0,1])}^\gamma
= C_{\gamma,T} ||u_0-v_0||_{L^1([0,1])}^\gamma.
\end{eqnarray*}
Point (iv) is checked similarly.
\end{preuve}

\section{Large time behavior}\label{ergo}

We now prove the uniqueness of the invariant measure.

\vip

\begin{preuve} {\it of Theorem \ref{uniinv}.}
Consider two invariant distributions $Q$ and $\tQ$ for (\ref{she}),
see Definition \ref{dinv}.
Let $u_0$ be $Q$-distributed and $\tu_0$ be $\tQ$-distributed.
Consider the corresponding (stationary) weak solutions $u,\tu$ to (\ref{she}).
Applying Theorem \ref{exstab}-(iv) and the Cauchy-Schwarz inequality, 
$\int_0^\infty K_sds <\infty$ a.s., where
\begin{eqnarray*}
K_s:=K(u(s),\tu(s))=
||b(u(s))-b(\tu(s))||_{L^1([0,1])} + 
||\sigma(u(s))-\sigma(\tu(s))||_{L^1([0,1])}^2.
\end{eqnarray*}
Using Lemma \ref{conc}, there is a sequence $(t_n)_{n\geq 1}$ such that
$K_{t_n}$ tends to $0$ in probability. Consider the function
$\phi(r)=r/(1+r)$ on $\rr_+$, and define
$\Psi:L^1([0,1])\times L^1([0,1]) \mapsto \rr_+$ as $\Psi(f,g)=\phi(K(f,g))$.
Then $\lim_n \E[\Psi(u(t_n),v(t_n))]=\lim_n \E[\phi(K_{t_n})]=0$.

We now apply Lemma \ref{cou}. The space $L^1([0,1])$ is Polish and
for each $n\geq 1$, $\cL(u(t_n))=Q$ and $\cL(\tu(t_n))=\tQ$. 
The function $\Psi$ is clearly continuous on $L^1([0,1])\times L^1([0,1])$,
(because $\sigma,b$ are Lipschitz-continuous).
Finally, $\Psi(f,g)>0$ for all $f\ne g$ (because 
$\Psi(f,g)=0$ implies that $b\circ f=b\circ g$ and 
$\sigma \circ f=\sigma \circ g$ 
a.e., whence $f=g$ a.e. since $(\sigma,b)$ is injective).
Lemma \ref{cou} thus yields $Q=\tQ$.
\end{preuve}

Finally, we give the

\vip

\begin{preuve} {\it of Theorem \ref{confl}.}
Point (ii) is immediately deduced from point (i). Let thus 
$u_0,v_0\in L^1([0,1])$
be fixed and let $u,v$ be the corresponding weak solutions
to (\ref{she}). We know from $(\cI)$, the Jensen inequality 
and Theorem \ref{exstab}-(iv) that a.s.,
\begin{eqnarray*}
\int_0^\infty \rho( ||u(t)-v(t)||_{L^1([0,1])}) dt &\leq&
\int_0^\infty || \rho(|u(t)-v(t)|)||_{L^1([0,1])} dt \\
&\leq&
\int_0^\infty \big|\big|  |b(u(t))-b(v(t))|+ |\sigma(u(t))-\sigma(v(t))|^2
\big|\big|_{L^1([0,1])} dt
<\infty.
\end{eqnarray*}
Using Lemma \ref{conc}, one may thus find an increasing
sequence $(t_n)_{n\geq 1}$
such that $\rho( ||u(t_n)-v(t_n)||_{L^1([0,1])})$ tends to $0$ in probability,
so that  $||u(t_n)-v(t_n)||_{L^1([0,1])}$ also tends to $0$ in probability
(because due to $\cI$, $\rho$ is strictly increasing and vanishes only at
$0$).
Next, we use Theorem \ref{exstab}-(iv) with e.g. $\gamma=1/2$ to get,
setting $\Delta_t=\sup_{[t,\infty)}||u(s)-v(s)||_{L^1([0,1])}$,
$$
\E\left[\left. \Delta_{t_n}^{1/2} \right\vert \cF_{t_n}\right] 
\leq C ||u(t_n)-v(t_n)||_{L^1([0,1])}^{1/2} \to 0 \hbox{ in probability}.
$$
We used here that conditionally on $\cF_{t_n}$, 
$(u(t_n+t,x))_{t\geq 0, x\in [0,1]}$ 
is a weak solution to (\ref{she}), starting from $u(t_n)$
(with a translated white noise).
Thus for any $\e>0$, using the Markov inequality
\begin{equation*}
P\left[\Delta_{t_n}>\e \right]
=  \E\left[ P\left( \left. \Delta_{t_n}>\e 
\right\vert \cF_{t_n} \right) \right] 
\leq \E\left[ \min\left(1, \e^{-1/2} 
\E\left[\left. \Delta_{t_n}^{1/2} \right\vert \cF_{t_n}\right]\right)\right],
\end{equation*}
which tends to $0$ as $n\to \infty$ by dominated convergence.
Consequently, as $n$ tends to infinity, 
\begin{equation}\label{cvtz}
\Delta_{t_n}
\hbox{ tends to $0$ in probability.}
\end{equation}
But a.s. $s\mapsto \Delta_s=\sup_{[s,\infty)} ||u(t)-v(t)||_{L^1([0,1])}$ is 
non-increasing,
and thus admits a limit as $s\to \infty$, which can be only $0$
due to (\ref{cvtz}).
\end{preuve}

\section{Toward the multi-dimensional case?}\label{dim2}

Consider now a bounded smooth domain $D \subset
\rr^d$, for some $d\geq 2$. Consider the (scalar) equation
\begin{equation}\label{sherd}
\partial_t u(t,x)=\Delta u(t,x) + b(u(t,x))+ \sigma(u(t,x))\dot W(t,x), \quad
t\geq 0,\; x\in D,
\end{equation}
with some Neumann boundary condition. 
Here $W(dt,dx)=\dot W (t,x)dtdx$ is a white noise 
on $[0,\infty)\times D$ based on $dtdx$.
We assume that
$\sigma,b:\rr \mapsto \rr$ are Lipschitz-continuous.

\vip

It is well known that the mild equation makes no sense in such a case,
since even if $\sigma(u)$ is bounded, $G_{t-s}(x,y)\sigma(u(s,y))$ does not
belong to $L^2([0,t]\times D)$. The existence of solutions
is thus still an open problem.
See however Walsh \cite{w} when $\sigma\equiv 1$,
$b(u)=\alpha u$
and Nualart-Rozovskii 
\cite{nr} when $\sigma(u)=u$, $b(u)=\alpha u$. In these works,
the authors manage to define some {\it ad-hoc} notion of solutions,
using that the equations can be solved more or less
explicitly. In the literature, one almost always
considers the simpler case where the noise $W$ is colored, see Da Prato-Zabczyk
\cite{dz}.

\vip

However the weak form makes sense: a predictable process
$u=(u(t,x))_{t\geq 0,x\in D}$ is a weak solution if a.s.,
\begin{equation}\label{crd} 
\hbox{for all } T>0, \quad
\sup_{[0,T]}||u(t)||_{L^1(D)}+ \int_0^T ||\sigma(u(t))||^2_{L^2(D)}dt <\infty
\end{equation}
and if for all function $\varphi\in C^2_b(D)$ 
(with Neumann conditions on $\partial D$), all $t\geq 0$, a.s.,
\begin{equation*}\label{weakrd}
\int_{D}u(t,x)\varphi(x)dx=\int_{D}u_0(x)\varphi(x)dx
+ \int_0^t\int_{D} [\{u(s,x)\Delta\varphi(x)+b(u(s,x))\}dxds+
\sigma(u(s,x))\varphi(x)W(ds,dx)].
\end{equation*}
Assume now that $\sigma(0)=b(0)=0$. Then $v\equiv 0$ is a 
weak solution. Furthermore, the estimate of Theorem \ref{exstab}-(iii)
{\it a priori} holds. Choosing $u_0\in L^1(D)$ and $v_0=0$, this would
imply (\ref{crd}). 
Unfortunately, we are not able to make this {\it a priori}
estimate rigorous. 

\vip

But following the proof of Proposition \ref{eql1}
and Corollary \ref{ineql1g}, one can easily check rigorously
the following result.
For $(e_k)_{k\geq 1}$  an orthonormal basis 
of $L^2(D)$, set $B^k_t= \int_0^t \int_{D} e_k(x) W(ds,dx)$. For 
$u_0 \in L^\infty(D)$ 
and $n\geq 1$, consider the solution (see Pardoux \cite{p}) to
$$
u^n(t,x)=u_0(x)+\intot [\partial_{xx}u^n(s,x) + b(u^n(s,x))]ds 
+ \sum_{k=1}^n \int_0^t \sigma(u^n(s,x))e_k(x)dB^k_s.
$$
Then if $\sigma(0)=b(0)=0$,
for any $\gamma\in (0,1)$, any $T>0$,
\begin{equation}\label{uin}
\E\left[ \sup_{[0,T]} ||u^n(t)||_{L^1(D)}^\gamma + \left\{
\int_0^T \sum_{k=1}^n \left( 
\int_{\rr^d} \sigma(u^n(t,x))e_k(x)dx \right)^2 ds \right\}^\gamma \right]
\leq C_{b,\gamma,T} ||u_0||_{L^1(D)}^\gamma,
\end{equation}
where the constant $C_{b,\gamma,T}$ depends only on $\gamma,T,b$
(the important fact is that it does not depend on $n$).
Passing to the limit formally in (\ref{uin}) would yield 
(\ref{crd}).
Unfortunately, (\ref{uin}) is not sufficient to ensure 
that the sequence $u^n$ is compact and tends, up to extraction of
a subsequence, to a 
weak solution $u$ to (\ref{sherd}). 
But this suggests
that, when $\sigma(0)=b(0)=0$, 
weak solutions to (\ref{sherd}) do exist and satisfy (\ref{crd}).

\section{Appendix}

First, we recall the following results on continuous local
martingales.

\begin{lem}\label{martingales}
Let $(M_t)_{t\geq 0}$ be a nonnegative continuous local martingale
starting from $m \in (0,\infty)$.
For all $\gamma \in (0,1)$, there exists a constant $C_\gamma$ (depending
only on $\gamma$) such that
$$
\E\left[\sup_{[0,\infty)} M_t^\gamma + 
\lb M \rb_\infty^{\gamma/2}  \right] \leq C_\gamma m^\gamma.
$$
\end{lem}

\begin{proof}
Classically (see e.g. Revuz-Yor \cite[Theorems 1.6 and 1.7 p 181-182]{ry} ), 
enlarging the probability space if necessary,
there is a 
standard Brownian motion $\beta$ such that $M_t=m+\beta_{\lb M \rb_t}$.
Denote now by $\tau_a=\inf\{t\geq 0;\; \beta_t=a\}$. Since $M$
is nonnegative, we deduce that
$$
\lb M\rb_\infty \leq \tau_{-m} \hbox{ and } \sup_{[0,\infty)} M_t
\leq m+\sup_{[0,\tau_{-m})} \beta_s.
$$ 
Thus we just have to prove that $\E[\tau_{-m}^{\gamma/2} ]
+\E[S_m^\gamma ] \leq C_\gamma m^\gamma$, where  $S_m=\sup_{[0,\tau_{-m})} \beta_s$.
\vip

First, for $x \geq 0$,  $P[S_m\geq x]=P[\tau_{x}\leq \tau_{-m}]=m/(m+x)$.
As a consequence, since $\gamma\in (0,1)$,
\begin{equation*}
\E[S_m^\gamma]= \int_0^\infty P[S_m^\gamma\geq x] dx
=\int_0^\infty \frac{m}{m+x^{1/\gamma}} dx
=m^\gamma \int_0^\infty \frac{1}{1+y^{1/\gamma}} dy = C_\gamma m^\gamma.
\end{equation*} 
Next, for $t\geq 0$, $P[\tau_{-m}\geq t]=P[\inf_{[0,t]}\beta_s > -m]$.
Recalling that $\inf_{[0,t]}\beta_s$ has the same law as $-\sqrt t |\beta_1|$, 
we get $P[\tau_{-m}\geq t]=P[|\beta_1|< m/\sqrt t ]$. Hence
\begin{equation*}
\E[\tau_{-m}^{\gamma/2}]= \int_0^\infty P [\tau_{-m}^{\gamma/2}\geq t ]dt=
\int_0^\infty P[|\beta_1|< m/ t^{1/\gamma}] dt
= \int_0^\infty P[(m/|\beta_1|)^\gamma>t ] dt = m^\gamma
\E\left[|\beta_1|^{-\gamma}\right].
\end{equation*}
This concludes the proof, since 
$\E\left[|\beta_1|^{-\gamma}\right]<\infty$ for $\gamma \in (0,1)$.
\end{proof}

Next, we state a technical result on a.s. converging integrals.

\begin{lem}\label{conc}
Let $(K_t)_{t\geq 0}$ be a nonnegative process. Assume that
$A_\infty= \int_0^\infty K_t dt <\infty$. Then one may find
a sequence $(t_n)_{n\geq 1}$ increasing to infinity such that
$K_{t_n}$ tends to $0$ in probability as $n\to \infty$.
\end{lem}

\begin{proof}
Consider a strictly increasing continuous concave function 
$\phi:\rr_+\mapsto [0,1]$ such that $\phi(0)=0$.
Using the Jensen inequality, we deduce that
$$
\frac{1}{T}\int_0^T \E[\phi(K_s)]ds= \E \left[ \frac{1}{T}\int_0^T \phi(K_s)ds
\right]\leq \E \left[\phi\left(\frac{1}{T}\int_0^T K_s ds\right)\right] 
\leq \E \left[\phi\left(\frac{A_\infty}{T}\right)\right],
$$
which tends to $0$ as $T\to \infty$ 
by the dominated convergence Theorem.
As a consequence, we may find a sequence $(t_n)_{n\geq 1}$ such that 
$\lim_n \E [\phi(K_{t_n})]=0$. The conclusion follows.
\end{proof}

Finally, we prove a technical result on coupling.

\begin{lem}\label{cou}
Consider two probability measures $\mu,\nu$ on a Polish space $\cX$.
Let $\Psi:\cX\times\cX \mapsto \rr_+$ be continuous and
assume that $\Psi(x,y)> 0$ for all $x\ne y$. If there is a sequence
of $\cX\times\cX$-valued random variables $(X_n,Y_n)_{n\geq 1}$
such that for all $n\geq 1$, $\cL(X_n)=\mu$ and $\cL(Y_n)=\nu$ and
if $\lim_n \E[\Psi(X_n,Y_n)]=0$, then $\mu=\nu$.
\end{lem}

\begin{proof}
The sequence of probability measures 
$(\cL(X_n,Y_n))_{n\geq 1}$ is obviously tight,
so up to extraction of a subsequence, we may assume that
$(X_n,Y_n)$ converges in law, to some $(X,Y)$. Of course,
$\cL(X)=\mu$ and $\cL(Y)=\nu$. Since $\Psi\land 1$ is continuous and bounded,
we deduce that $\E[\Psi(X,Y)\land 1]=\lim_n \E[\Psi(X_n,Y_n)\land 1]=0$, 
whence $\Psi(X,Y)=0$ a.s. By 
assumption, this implies that $X=Y$ a.s., so that $\mu=\nu$.
\end{proof}


\begin{thebibliography}{99}
\rm

\bibitem{bgp}{V. Bally, I. Gyongy, E. Pardoux, {\it White noise driven 
parabolic SPDEs with measurable drift,} 
J. Funct. Anal. 120, no. 2, 484--510, 1994. }

\bibitem{bp}{V. Bally, E. Pardoux, {\it Malliavin calculus for white noise 
driven parabolic SPDEs,}  Potential Anal.  9,  no. 1, 27--64, 1998.}

\bibitem{bms}{V. Bally, A. Millet, M. Sanz-Sol\'e, {\it Approximation and 
support in H\"older norm for parabolic stochastic partial differential
equations}, Annals of Probability 23, no 1, 178--222, 1995.}

\bibitem{c}{S. Cerrai, {\it Stabilization by noise for a class
of stochastic reastion-diffusion equations}, Probab. Theory Relat. Fields
133, 190-214, 2005.}

\bibitem{dz}{G. Da Prato, J. Zabczyk, {\it Stochastic equations in infinite 
dimensions,}, Cambridge University Press, 1992.}

\bibitem{dp}{C. Donati-Martin, E. Pardoux, {\it White nise driven SPDEs
with reflection}, Probab. Theory Relat. Fields 95, 1-24, 1993.}

\bibitem{gg}{D. Gatarek, B. Goldys, {\it On weak solutions of stochastic 
equations in Hilbert spaces}, Stochastics Stochastics Rep.  46,  no. 1-2, 
41--51, 1994.}

\bibitem{nr}{D. Nualart, B. Rozovskii, {\it Weighted stochastic Sobolev
spaces and bilinear SPDEs driven by space-time white noise}, J. Fucnt.
Anal. 149, 200-225, 1997.}

\bibitem{m}{C. Mueller, {\it Coupling and invariant measures for the heat 
equation with noise}, Ann. Probab., 21, 2189-2199, 1993.}

\bibitem{p}{E. Pardoux, {\it Stochastic partial differential equations
and filtering of diffusion processes}, Stochastics 3, 127-167, 1979.}

\bibitem{ry}{D. Revuz, M. Yor, {\it Continuous martingales and Brownian
motion,} third edition, Springer Verlag, 1999.}

\bibitem{s}{R. Sowers, {\it Large deviations for the invariant measure
of a reaction-diffusion equation with non-Gaussian perturbations},
Probab. Theory Relat. Fields 92, 393-421, 1992.}

\bibitem{w}{J.B. Walsh, {\it An introduction to stochastic partial 
differential equations,}  Ecole d'\'et\'e de probabilit\'es de 
Saint-Flour, XIV---1984,  265--439, Lecture Notes in Math., 1180, 
Springer, Berlin, 1986.}

\end{thebibliography}
\end{document}